\documentclass[review]{article}
\usepackage{amsmath,amsfonts,amssymb,amsthm,amsxtra}
\usepackage{tipa}
\usepackage[T1]{fontenc}
\usepackage[latin1,latin9,ansinew]{inputenc}
\usepackage{graphicx,latexsym}  
\usepackage[a4paper,top=3cm,left=3cm,right=3cm]{geometry}
\usepackage[capitalise]{cleveref}

\def \C{\mathbb{C}}
\def \K{\mathbb{K}}
\def \LL{\mathbb{L}}
\def \FF{\mathbb{F}}
\def \N{\mathbb{N}}

\def \Q{\mathbb{Q}}
\def \R{\mathbb{R}}

\newtheorem{thm}{Theorem}

\newtheorem{lemma}{Lemma}
\newtheorem{prop}{Proposition}

\newtheorem{Rem}{Remark}

\def \a{{\"a}}

\usepackage{lipsum}

\usepackage{float}
\providecommand{\keywords}[1]
{
  \small	
  \textbf{\textit{Keywords:\,}} #1
}

\title{Abstract, keywords and references template}
\author{Gessica Alecci$^{1}$, Carsten Elsner$^{2}$  \\
        \small $^{1}$DISMA, Politecnico di Torino \\
        \small $^{2}$Fachhochschule fur die Wirtschaft, Hannover \\
}

\title{On a criterion for algebraic independence and its variants}
\date{}

\begin{document}
\maketitle

\begin{abstract}
From around 2010 onward, Elsner et al.\,developed and applied a method in which the algebraic independence of $n$ quantities $x_1,\dots,x_n$ over a field is transferred to further $n$ quantities $y_1,\dots,y_n$ by means of a system of polynomials in $2n$ variables $X_1,\dots,X_n,$ $Y_1,\dots,Y_n$. In this paper, we systematically study and explain this criterion and its variants. 
\end{abstract}


\section{Introduction} \label{sec:1}
The transcendence of the circular number $\pi$ and that of the number $e=\exp(1)$ has been known since the end of the 19th century but the question of the algebraic independence of these two numbers over ${\Q}$ has still not been answered. It concerns the exclusion of the existence of a non-identical vanishing polynomial
$P(X,Y)$ with rational coefficients such that $P(\pi,e)=0$.

The theorem of Lindemann - Weierstrass (1885), from which
the transcendence of $\pi$ and of $e$ can be derived, is the beginning of a general theory on algebraic independence of complex numbers over ${\Q}$. In one of its equivalent formulations
this theorem states that in the case of the {\em linear\/} independence of algebraic numbers
$\alpha_1,\dots,\alpha_n$ over ${\Q}$, the numbers $e^{\alpha_1},\dots,e^{\alpha_n}$ 
are {\em algebraically\/} independent over ${\Q}$ \cite{Shidlovskii}.

An additional significant achievement is the theorem of Gelfond-Schneider that states the transcendence of $\alpha^{\beta}$ when $\alpha$ and $\beta$ are 
algebraic over ${\Q}$, assuming that $\alpha \not= 0,1$ and $\beta \not\in {\Q}$ \cite{Shidlovskii}. Another important result is the Baker's Theorem on linear forms of logarithms that states that given $\alpha_1,\dots,\alpha_n$ algebraic numbers different from zero such that $\log \alpha_1,\dots,\log \alpha_n$ are linearly independent over the rational numbers, then the numbers $1,\log \alpha_1,\dots,\log \alpha_n$
are linearly independent over the field of all algebraic numbers \cite{Baker}.

In 1916, S.\,Ramanujan \cite{Ramanujan} defined the series
\[S_{j+1}(x) \,:=\, \frac{\zeta(-2j-1)}{2} + \sum_{n=1}^{\infty} \frac{n^{2j+1}x^n}{1-x^n}\,\]
where $\zeta(s)$ is the Riemann zeta function. Let
\begin{equation}
P(x) \,:=\, -24S_1(x) \,,\qquad Q(x) \,:=\, 240S_3(x) \,,\qquad R(x) \,:=\, -504S_5(x) \,.
\label{RamaFunctions}
\end{equation}
In 1996, Y.\,Nesterenko \cite{Nesterenko1} proved that for every complex number $x$ with $0<|x|<1$, the set 
\[\big\{  x,P(x),Q(x),R(x)\big\} \] 
contains at least three numbers that are algebraically independent over ${\Q}$.

In the years since around 2010, the second named author, in collaboration with the Japanese mathematicians I.\,Shiokawa and Sh.\,Shimomura, obtained numerous results on algebraic independence or dependence of number sets. They started from a set of known algebraic independent numbers and they obtained a second set of numbers, where the numbers in both sets satisfy a system of polynomial equations. Then they developed a criterion, which allows to decide on the algebraic independence of the numbers in the second set. Different variants of this criterion have so far appeared in various publications as essential tools \cite{Elsner1}, \cite{Elsner3} and \cite{Elsner8}. However, they have not yet been systematically summarized and explained in a journal, taking into account their interrelationships. The goals of this paper is to systematically explain this criterion and its variants including many examples, obtained by the second name author and his collaborators, which can be found in the bibliography. Section \ref{sec:2} contains the statements of the theorems regarding this criterion, namely Theorems \ref{Thm1} and \ref{Thm2}, whereas Theorems \ref{Thm3} and \ref{Thm4} contain sufficient conditions for the independence criterion.
The theorems are proven in Section \ref{sec:3}. Finally, in Section \ref{sec:4} we present some past results obtained by applications of the variants of the algebraic independence criterion. Many of these results are going back to the above-mentioned algebraic independence result for $P(x),Q(x)$ and $R(x)$ over ${\Q}$ for a non-vanishing algebraic number $x$ with $|x|<1$, obtained by Y.\,Nesterenko. Other results rely on the Lindemann-Weierstrass Theorem.

\section{Statement of the criterion and its variants} \label{sec:2}
We summarize the different variants of the criterion taking into account their interrelationships. 
\begin{thm}
Let ${\K}$ be a field with ${\Q} \subseteq {\K} \subseteq {\C}$. Furthermore, it is assumed that the numbers
$x_1,\dots ,x_n \in {\C}$ and $y_1,\dots ,y_n \in {\C}$ satisfy a system
\begin{equation}
f_j(x_1,\dots ,x_n,y_1,\dots ,y_n) \,=\, 0 \qquad (j=1,\dots ,n)
\label{10}
\end{equation}
of equations with polynomials $f_j(X_1,\dots ,X_n,Y_1,\dots ,Y_n) \in {\K}[X_1,\dots ,X_n,$ $Y_1,\dots ,Y_n]$ for $j=1,
\dots ,n$. If the numbers $x_1,\dots ,x_n$ are algebraically independent over ${\K}$ and 
\begin{equation}
\mathrm{det}_n \,\Big( \,\frac{\partial f_j}{\partial X_i} (x_1,\dots ,x_n,y_1,\dots ,y_n)\,\Big) \,\not= \,0 
\label{20}
\end{equation}
holds, then the numbers $y_1,\dots ,y_n$ are algebraically independent over ${\K}$.
\label{Thm1}
\end{thm}
\begin{Rem}
 We interpret $j$ as the row number and $i$ as the column number. From now on, we denote by capital letters the variables and by lowercase letters the numerical values they assume. Moreover, by det$_k$ we mean the determinant of a $k \times k$ matrix.
\end{Rem}\label{remarkindex}
\begin{thm}
Let ${\K}$ be a field with ${\Q} \subseteq {\K} \subseteq {\C}$. Furthermore, it is assumed that the numbers
$x_1,\dots ,x_n \in {\C}$ and $y_1,\dots ,y_n \in {\C}$ satisfy a system
\begin{equation}
y_j \,=\, T_j(x_1,\dots ,x_n) \qquad (j=1,\dots ,n)
\label{30}
\end{equation}
of equations with polynomials $T_j(X_1,\dots ,X_n) \in {\K}[X_1,\dots ,X_n]$ for $j=1,\dots ,n$. If the
numbers $x_1,\dots ,x_n$ are algebraically independent over ${\K}$ and 
\begin{equation}
\mathrm{det}_n \,\Big( \,\frac{\partial T_j}{\partial X_i} (x_1,\dots ,x_n)\,\Big) \,\not= \,0 
\label{40}
\end{equation}
holds, then the numbers $y_1,\dots ,y_n$ are algebraically independent over ${\K}$.
\label{Thm2}
\end{thm}
\begin{Rem}
 Under the conditions of Theorem~\ref{Thm2}, the non-vanishing of the determinant in (\ref{40}) is not only sufficient but 
also necessary for the algebraic independence of $y_1,\dots,y_n$ over ${\K}$, but there are no interesting application of it.
\end{Rem}\label{remarknecessarycond}
In numerous applications, each of the numbers $y_1,\dots,y_n$ can be represented as the value of a {\em rational
function\/} at the position $(x_1,\dots,x_n)$:
\begin{equation}
y_j \,=\, \frac{T_j(x_1,\dots,x_n)}{U_j(x_1,\dots,x_n)} \qquad (j=1,\dots,n) \,.
\label{100}
\end{equation}
Here $T_j$ and $U_j$ are non-identical vanishing polynomials from the ring ${\K}[X_1,\dots,X_n]$. The polynomials
\begin{equation}
f_j(X_1,\dots,X_n,Y_j) \,:=\, Y_jU_j(X_1,\dots,X_n)-T_j(X_1,\dots,X_n)
\label{110}
\end{equation}
for $j=1,\dots,n$ thus vanish at the point $(x_1,\dots,x_n,y_j)$. 
\begin{thm}
Let ${\K}$ be a field with ${\Q} \subseteq {\K} \subseteq {\C}$. Furthermore, we assume that the
numbers $x_1,\dots ,x_n \in {\C}$ and $y_1,\dots ,y_n \in {\C}$ satisfy a system
\begin{equation*}
y_j \,=\, R_j(x_1,\dots ,x_n) \qquad (j=1,\dots ,n)
\label{120}
\end{equation*}
of equations with rational functions $R_j(X_1,\dots ,X_n)$ in the extension field $ {\K}(X_1,\dots ,X_n)$ for $j=1,\dots ,n$. If the
numbers $x_1,\dots ,x_n$ are algebraically independent over ${\K}$ and 
\begin{equation}
\mathrm{det}_n \,\Big( \,\frac{\partial R_j}{\partial X_i} (x_1,\dots ,x_n)\,\Big) \,\not= \,0 
\label{130}
\end{equation}
holds, then the numbers $y_1,\dots ,y_n$ are algebraically independent over ${\K}$.
\label{Thm3}
\end{thm}
\begin{thm}
Let ${\K}$ be a field with ${\Q} \subseteq {\K} \subseteq {\C}$. Let $m$ and $n$ be positive integers with
$1\leq m<n$. Furthermore, we assume that the numbers $x_1,\dots ,x_n \in {\C}$ and $y_1,\dots ,y_m \in
{\C}$ satisfy a system
\begin{equation}
f_j(x_1,\dots ,x_n,y_1,\dots ,y_m) \,=\, 0 \qquad (j=1,\dots ,m)
\label{150}
\end{equation}
of equations with polynomials $f_j(X_1,\dots ,X_n,Y_1,\dots ,Y_m) \in {\K}[X_1,\dots ,X_n,$ $Y_1,\dots ,Y_m]$ for
$j=1,\dots ,m$. If the numbers $x_1,\dots ,x_n$ are algebraically independent over ${\K}$ and if
\begin{equation}
\mathrm{det}_m \,\Big( \,\frac{\partial f_j}{\partial X_i} (x_1,\dots ,x_n,y_1,\dots ,y_m)\,\Big) \,\not= \,0 
\label{160}
\end{equation}
holds with $1\leqslant i,j\leqslant m$, then the numbers $y_1,\dots ,y_m$ are algebraically independent over the field ${\K}(x_{m+1},\dots,x_n)$.
\label{Thm4}
\end{thm}

\section{Proofs of the algebraic criterion and its variants} \label{sec:3}
There exist a lot of different proofs of Theorem \ref{Thm1} that make use of the different tools. A first proof, that will be shown below, makes use of commutative algebra and a lemma on separable algebraic field extensions. Other proofs make use respectively of the 
projection theorem of Tarski - Seidenberg and the concept of semi-algebraic sets (in the case of a real field ${\K}$), differential forms or isomorphism of fields, see \cite[Example\,1]{Elsner2}. Finally, one last proof is an elementary one, based on the essential idea to embed a given zero of a polynomial with at least two variables into a path on which the polynomial  vanishes identically. In particular, the differential forms allow us to prove the converse statement of Theorem\,\ref{Thm2}. We omit the proof in this paper, because this result is only of theoretical interest and has no interesting application.

\begin{prop}
Let ${\LL}$ be a field with ${\Q} \subseteq {\LL} \subseteq {\R}$. Furthermore, let the point $(x_1,\dots ,x_n) \in {\R}^n$ be an isolated zero of the system of
equations
\begin{equation*}
P_j(X_1,\dots,X_n) \,=\, 0 \qquad (j=1,\dots ,n)\,,
\label{50}
\end{equation*}
which is formed with polynomials $P_j(X_1,\dots ,X_n) \in {\LL}[X_1,\dots ,X_n]$ for $j=1,\dots,n$. Then the numbers $x_1,\dots ,x_n$ are algebraic over ${\LL}$.
\label{Prop1}
\end{prop} 
The following lemma, which will replace Proposition\,\ref{Prop1}, is essential in the proof of Theorem\,\ref{Thm1} that we show in this paper.
\begin{lemma}
Let ${\LL}$ denote a field, ${\FF}$ a field extension of ${\LL}$. We assume that the quantities $x_1,\dots,x_n \in \FF$ satisfy a system
\[ P_j(x_1,\dots ,x_n) \,=\, 0 \qquad (j=1,\dots ,n) \]
of equations with polynomials $P_j(X_1,\dots,X_n) \in {\LL}[X_1,\dots,X_n]$. If
\[\mathrm{det}_n\,\Big( \,\frac{\partial P_j}{\partial X_i} (x_1,\dots ,x_n)\,\Big) \,\not= \,0 \]
holds, then ${\LL}(x_1,\dots,x_n)$ is a separable algebraic field extension of ${\LL}$.
\label{Lem3}
\end{lemma}
{\em Proof.\/}
This lemma is the Corollary of Theorem~40 in \cite[p.\,126]{Zariski}. It also results from Proposition~5.3 in \cite[p.\,371]{Lang}. \hfill \qed \vspace*{5pt} \\
{\em Proof of Theorem\,\ref{Thm1}.\/}\,
Using the polynomials $f_1,\dots ,f_n$ from Theorem~\ref{Thm1}, we define new polynomials $P_j$ for $j=1,\dots,n$ by
\begin{equation}
P_j(X_1,\dots ,X_n) \,:=\, f_j(X_1,\dots ,X_n,y_1,\dots ,y_n) \in {\K}(y_1,\dots ,y_n)[X_1,\dots ,X_n] \,.
\label{70}
\end{equation}
They are constructed using the variables $X_1,\dots ,X_n$, and they all vanish at the point $(x_1,\dots,$ $x_n)$ because of (\ref{10}).
The determinant condition (\ref{20}) takes the form
\begin{equation}
\mathrm{det}_n \,\Big( \,\frac{\partial P_j}{\partial X_i} (x_1,\dots ,x_n)\,\Big) \,\not= \,0 \,.
\label{80}
\end{equation}
With $\FF={\C}$, ${\LL}={\K}(y_1,\dots,y_n)$ and the setting of polynomials $P_j(X_1,\dots,X_n)$ $(j=1,\dots,n)$ as in (\ref{70}),  all the conditions of Lemma~\ref{Lem3} 
including (\ref{80}) are fulfilled. So, by Lemma\,\ref{Lem3}, ${\LL}(x_1,\dots ,x_n)$ is an algebraic field extension over ${\LL}$. \hfill \vspace*{5pt} \\
We denote by $tr.deg\big( {\K}_2 : {\K}_1 \big)$ the transcendence degree of the field extension of ${\K}_2$ over ${\K}_1$. Then we have 
\[tr.deg\big( {\LL}(x_1,\dots ,x_n) : {\LL} \big) \,=\, 0 \quad \mbox{and} \quad tr.deg \big( {\K}(x_1,\dots ,x_n) : {\K} \big) \,=\,
n \,; \]
the latter due to the condition on the algebraic independence of $x_1,\dots ,x_n$ over ${\K}$ in Theorem~\ref{Thm1}. Trivially, ${\K} \subseteq {\LL}$ implies that 
$tr.deg ({\LL}(x_1,\dots ,x_n):{\K}) \geq n$. If the fields ${\K}_1$, ${\K}_2$ and ${\K}_3$ form a field tower ${\K}_1 \subseteq {\K}_2 \subseteq {\K}_3$,
we know according to the chain theorem for degrees of transcendence that we have the identity
\begin{equation}
tr.deg\big( {\K}_3 : {\K}_1 \big) \,=\, tr.deg\big( {\K}_3 : {\K}_2 \big) + tr.deg\big( {\K}_2 : {\K}_1 \big) \,.
\label{90}
\end{equation}
Applying this equation to the field tower ${\K} \subseteq {\LL}={\K}(y_1,\dots ,y_n) \subseteq
{\LL}(x_1,\dots ,x_n)$, we obtain
\begin{eqnarray*}
n &\,\leq \,& tr.deg\big( {\LL}(x_1,\dots ,x_n) : {\K} \big) \,=\, tr.deg\big( {\LL}(x_1,\dots ,x_n) : {\LL} \big) +
tr.deg\big( {\LL} : {\K} \big) \\
&\,=\,& tr.deg\big( {\LL} : {\K} \big) \,=\, tr.deg\big( {\K}(y_1,\dots ,y_n) : {\K} \big) \,\leq \, n \,.
\end{eqnarray*}
By $tr.deg\big( {\K}(y_1,\dots ,y_n) : {\K} \big) =n$, Theorem~\ref{Thm1} is proved. \hfill \qed \vspace*{8pt} \\

{\em Proof of Theorem\,\ref{Thm2}.\/}\, Setting
\[f_j\big( X_1,\dots,X_n,Y_j \big) \,:=\, T_j\big( X_1,\dots,X_n \big) - Y_j \qquad (1\leqslant j \leqslant n)\,,\]
Theorem\,\ref{Thm2} follows immediately from Theorem\,\ref{Thm1}. \hfill \qed \vspace*{5pt} \\

{\em Proof of Theorem\,\ref{Thm3}.\/}\,
We denote the determinant in (\ref{130}) by $\Delta_1$ and assume first that $\Delta_1 \not=0$. 
Furthermore, we express the rational functions $R_j$ by
\[R_j(X_1,\dots,X_n) \,=\, \frac{T_j(X_1,\dots,X_n)}{U_j(X_1,\dots,X_n)} \qquad (j=1,\dots,n)\,, \]
using the polynomials $T_j$ and $U_j$ already introduced in (\ref{100}). The algebraic independence of $y_1,\dots,y_n$ over ${\K}$ can be
 proven with Theorem~\ref{Thm1}, where the non-vanishing of the Jacobi determinant  must be shown for the functions $f_j$ from (\ref{110}):
\begin{eqnarray*}
\Delta_2 &\,:=\,& \mathrm{det}_n \,\Big( \,\frac{\partial f_j}{\partial X_i} \,\Big) (x_1,\dots ,x_n,y_1,\dots ,y_n) \nonumber \\
&\,=\,& \mathrm{det}_n \,\Big( \,y_j\frac{\partial U_j}{\partial X_i} - \frac{\partial T_j}{\partial X_i}\,\Big)
(x_1,\dots ,x_n) \nonumber \\
&\,=\,& \mathrm{det}_n \,\Big( \,\frac{T_j}{U_j}\frac{\partial U_j}{\partial X_i} - \frac{\partial T_j}{\partial X_i}\,\Big)
(x_1,\dots ,x_n) \nonumber \\
&\,=\,& {(-1)}^nU_1\cdots U_n \mathrm{det}_n \,\left( \,\frac{1}{U_j^2}\Big( \,\frac{\partial T_j}{\partial X_i} U_j -
T_j\frac{\partial U_j}{\partial X_i} \Big) \,\right)(x_1,\dots ,x_n) \nonumber \\
&\,=\,& {(-1)}^nU_1\cdots U_n \mathrm{det}_n \,\Big( \,\frac{\partial R_j}{\partial X_i}\,\Big) (x_1,\dots ,x_n) \nonumber \\
&\,=\,& {(-1)}^nU_1\cdots U_n (x_1,\dots ,x_n) \Delta_1 \,.
\label{140}
\end{eqnarray*}
From $\Delta_1 \not= 0$ we have $\Delta_2 \not= 0$. The $x_1,\dots,x_n$ are considered as algebraically
independent over ${\K}$ so that the polynomial $Q_1\cdots Q_n$ does not vanish identically. Then,
$U_1\cdots U_n (x_1,\dots ,x_n) \not= 0$ is guaranteed. Theorem~\ref{Thm1} thus proves the
algebraic independence of $y_1,\dots,y_n$ over ${\K}$.  \hfill \qed \vspace*{8pt} \\
{\em Proof of Theorem\,\ref{Thm4}.\/}\,
In addition to the equations in (\ref{150}) we introduce the polynomials
\[f_j(X_j,Y_j) \,:=\, X_j-Y_j \,\in \, {\K}[X_j,Y_j] \qquad (j=m+1,\dots,n) \]
so that for $y_j:=x_j$ $(j=m+1,\dots,n)$ the polynomial $f_j$ vanishes at the point
$(x_j,y_j)$. Theorem~\ref{Thm1} can now be applied for the proof of the algebraic independence of 
$y_1,\dots,y_n$ over ${\K}$. We compute the Jacobi determinant of the functions $f_1,\dots,f_n$ at
the position $(x_1,\dots,x_n,y_1,\dots,y_n)$ and apply the condition from (\ref{160}):
\begin{eqnarray*}
&& \mathrm{det}_n \,\Big( \,\frac{\partial f_j}{\partial X_i} (x_1,\dots ,x_n,y_1,\dots ,y_n)\,\Big) \\
&\,=\,& \mathrm{det}_n \,\left( \,\begin{array}{cccccc}
\displaystyle \frac{\partial f_1}{\partial X_1} & \cdots & \displaystyle \frac{\partial f_1}{\partial X_m} &
\displaystyle \frac{\partial f_1}{\partial X_{m+1}} & \cdots & \displaystyle \frac{\partial f_1}{\partial X_n} \\
\vdots & \cdots & \vdots & \vdots & \cdots & \vdots \\
\displaystyle \frac{\partial f_m}{\partial X_1} & \cdots & \displaystyle \frac{\partial f_m}{\partial X_m} &
\displaystyle \frac{\partial f_m}{\partial X_{m+1}} & \cdots & \displaystyle \frac{\partial f_m}{\partial X_n} \\
0 & \cdots & 0 & 1 & \cdots & 0 \\
\vdots & \cdots & \vdots & \cdots & \cdots & \vdots \\
0 & \cdots & 0 & 0 & \cdots & 1
\end{array} \,\right) (x_1,\dots,x_n,y_1,\dots,y_m) \\
&\,=\,& \mathrm{det}_n \,\left( \,\begin{array}{ccc}
\displaystyle \frac{\partial f_1}{\partial X_1} & \cdots & \displaystyle \frac{\partial f_1}{\partial X_m} \\
\vdots & \cdots & \vdots \\
\displaystyle \frac{\partial f_m}{\partial X_1} & \cdots & \displaystyle \frac{\partial f_m}{\partial X_m}
\end{array}
\,\right) (x_1,\dots,x_n,y_1,\dots,y_m) \quad \not= \, 0 \,.
\end{eqnarray*}
Thus, according to Theorem~\ref{Thm1} and to $y_j=x_j$ for $j=m+1,\dots,n$, the algebraic independence of the numbers
$y_1,\dots,y_m,x_{m+1},\dots,x_n$ over ${\K}$ is proven. This implies 
\[ tr.deg \big( {\K}(y_1,\dots,y_m,x_{m+1},\dots,x_n):{\K} \big)  \,=\,n\,.\] 
Taking into account the presupposed algebraic independence of $x_1,\dots,x_n$ over ${\K}$,
we know that the equation $tr.deg \big( {\K}(x_{m+1},\dots,x_n):{\K} \big) =n-m$ holds. Again we make use of the
chain rule from (\ref{90}), applied to the fields ${\K} \,\subseteq \, {\K}(x_{m+1},\dots,x_n) \,\subseteq \, 
{\K}(y_1,\dots,y_m,x_{m+1},\dots,x_n)$. Thus we obtain
\[tr.deg \big( {\K}(y_1,\dots,y_m,x_{m+1},\dots,x_n):{\K}(x_{m+1},\dots,x_n) \big) \,=\, m \,,\]
which is the statement of Theorem~\ref{Thm4}. \hfill \qed 

\section{Past results obtained with the algebraic independence criterion and its variants} \label{sec:4}
We summarize some results that have been obtained with the algebraic independence criterion and its variants. It is only a small selection of scattered published results, but it is intended 
to demonstrate the spectrum of applications of the method.
\begin{itemize}
\item[1.)] In the following we state some applications of Theorem\, \ref{Thm1}.
\begin{itemize}
\item Let $j\geqslant 0$ be an integer. There are 16 families of $q$-series like the three families 
\begin{eqnarray}
A_{2j+1}(q) &=& \sum_{n=1}^{\infty} \frac{n^{2j+1}q^{2n}}{1-q^{2n}}\,,
\label{Ref10} \\
B_{2j+1}(q) &=&= \sum_{n=1}^{\infty} \frac{{(-1)}^{n+1}n^{2j+1}q^{2n}}{1-q^{2n}} \,,
\nonumber \\
C_{2j+1}(q) &=& \sum_{n=1}^{\infty} \frac{n^{2j+1}q^n}{1-q^{2n}}\,,
\nonumber
\end{eqnarray}
which are generated from the Fourier expansion of the Jacobian elliptic functions. In \cite[Theorem\,2]{Elsner1} sets of three such algebraically independent 
$q$-series for an algebraic number $q$ with $0<|q|<1$ are characterized. For example, $A_{2j+1}(q),B_{2j+1}(q)$ and $C_{2j+1}(q)$ are algebraically independent 
over ${\Q}$.
\item The paper \cite{Elsner6} was the starting point for a series of papers on algebraic dependence and independence results of the values of theta-constants. In
Theorem\,1.1 of this paper it is shown among other things that for any algebraic number $q$ with $0<|q|<1$ the two numbers $\theta(q)$ and $\theta(q^{2^m})$ with 
an integer $m\geqslant 1$ are algebraically independent over ${\Q}$, where
\[\theta(q) \,:=\, 1+2\sum_{\nu =1}^{\infty} q^{\nu^2} \,.\]  
\item Let $q\geq 3$ be an integer and consider the power series
\[e_r(z) \,:=\, \sum_{\scriptsize \begin{array}{c} n=0 \\ n \equiv r \pmod q \end{array}}^{\infty} \frac{z^n}{n!} \qquad (r=0,1,\dots,q-1)\,.\] 
In \cite[Theorem 1]{Elsner8} it is shown that for any nonzero algebraic number $\alpha$, among $q$ numbers $e_0(\alpha),e_1(\alpha),\dots,e_{q-1}(\alpha)$ any
$\varphi(q)$ are algebraically independent over ${\Q}$, where $\varphi$ denotes Euler's totient. 
\end{itemize}
\item[2.)] Here there are some applications of Theorem\,\ref{Thm2}.
\begin{itemize}
\item  Let $q$ be an algebraic number with $0<|q|<1$. Then the three numbers $A_1(q)$, $A_{2i+1}(q)$ and $A_{2j+1}(q)$
from the series in (\ref{Ref10}) with $1\leqslant i<j$ and $(i,j)\not= (1,3)$ are algebraically independent over ${\Q}$. For $A_1(q)$, $A_3(q)$ and $A_7(q)$ there
is the algebraic relation $A_7(q)=A_3(q) + 120A_3^2(q)$ \cite[Theorem\,2]{Elsner2}.
\item  Let $F_n$ denote the Fibonacci numbers with $F_0=0$, $F_1=1$ and $F_{n}=F_{n-1}+F_{n-2}$ for all $n\geq 2$. Let $s_1,s_2,s_3$ be distinct positive integers. 
It follows from \cite[Theorem\,1.1]{Elsner3} that the three numbers
\begin{equation}
\zeta_{Fib}(2s_i) \,:=\,\sum_{n=1}^{\infty}\frac{1}{F_n^{2s_i}} \qquad (i=1,2,3) 
\label{Ref20}
\end{equation}
are algebraically independent over ${\Q}$ if and only if at least one of $s_1,s_2,s_3$ is even.
\item In the PhD thesis \cite{Stein} the previous result on the series given in (\ref{Ref20}) is generalized to reciprocal sums of sequences of integers satisfying 
a linear three-term recurrence formula,
see \cite[Theorem\,5.3]{Stein}. 
\end{itemize}
\item[3.)] Finally, we sit some applications of Theorem\,\ref{Thm3}.
\begin{itemize}
\item The Ramanujan function $P(x)$ is given by $P(q^2)=1-24A_1(q)$ with $A_1(q)$ defined in (\ref{Ref10}), see also (\ref{RamaFunctions}). By 
\cite[Theorems\,5.1]{Elsner5} we know that for every algebraic number $q$ with $0<|q|<1$ any three numbers in the set
\[\big\{ \,P(q),\,P(q^2),\,P(q^5),\,P(q^{10})\,\big\} \]
are algebraically independent over ${\Q}$, and the four numbers are not. By Theorem\,6.1 in the same paper \cite{Elsner5}, for positive integers $a$ and $b$ with
$a^2-4b<0$, any two of the numbers
\[\sum_{n=0}^{\infty} \frac{1}{n^2+acn+bc^2} \qquad (c=1,2,3,\dots) \]
are algebraically independent over ${\Q}$ and any three of them are not. \vspace*{5pt} \\
The Lucas numbers $L_n$ are defined recursively by $L_0=2$, $L_1=1$ and $L_{n}=L_{n-1}+L_{n-2}$ for all $n\geq 2$. Let $f_s(\alpha)$ and $g_s(\alpha)$ be power
series defined by
\[f_s(z) \,=\, \sum_{n=0}^{\infty} F_n^s \frac{z^n}{n!} \quad \mbox{and} \quad g_s(z) \,=\, \sum_{n=0}^{\infty} L_n^s \frac{z^n}{n!} \,,\]
where $s\in {\N}$. Let $\alpha$ be a nonzero algebraic number. Then, by \cite[Theorem\,5]{Elsner8}, all the numbers in the set $\{ f_s(\alpha):s\in {\N} \} \cup
\{ g_s(\alpha):s\in {\N} \}$ are distinct and any two are algebraically independent over ${\Q}$. Moreover, any three functions in the set $\{ f_s(z):s\in {\N}
\} \cup \{ g_s(z):s\in {\N} \}$ are algebraically dependent over ${\Q}$. A similar statement holds by Theorem\,6 in the same paper \cite{Elsner8} for two power
series defined by
\[f_{a,b}(z) \,=\, \sum_{n=0}^{\infty} F_{an+b} \frac{z^n}{n!} \quad \mbox{and} \quad g_{a,b}(z) \,=\, \sum_{n=0}^{\infty} L_{an+b}^s \frac{z^n}{n!} \,.\]
If $\alpha$ is a nonzero algebraic number, then any two numbers in the set $\{ f_{a,b}(\alpha):a\in {\N},b\in {\N}_0 \}$ are algebraically independent over ${\Q}$.
Moreover, any three functions in the set $\{ f_{a,b}(z):s\in {\N} \}$ are algebraically dependent over ${\Q}$. The same statements hold for the power series 
$g_{a,b}(z)$.
\end{itemize}
\item[4.)] Regarding Theorem\,\ref{Thm4}, we show the application of it by an example. \vspace*{5pt} \\
We investigate the values of the two series 
\begin{eqnarray*}
y_1 &:=& \zeta_{Fib}(4) \,=\, \sum_{n=1}^{\infty} \frac{1}{F_n^4} \,=\, 2.076730850 \dots \,,\\
y_2 &:=& \zeta_{Fib}(8) \,=\, \sum_{n=1}^{\infty} \frac{1}{F_n^8} \,=\, 2.004061286 \dots \,.
\end{eqnarray*}
One can express $y_1$ and $y_2$ each by three parameters which are associated with the complete elliptic
integrals of the first and second kind,
\[K(k) \,:=\, \int_0^1 \frac{dt}{\sqrt{(1-t^2)(1-k^2t^2)}} \,,\quad
E(k) \,:=\, \int_0^1 \sqrt{\frac{1-k^2t^2}{1-t^2}}\,dt \,,\]
and with the underlying modulus $k$. This modulus is given by the uniquely determined real number $k$
in the interval $[0,1]$ such that
\[-\frac{2}{\pi} \log \Big( \,\frac{\sqrt{5}-1}{2}\, \Big) \,=\, \frac{K\big( \sqrt{1-k^2} \,\big)}{K(k)} \,.\]
Moreover, we introduce three parameters $x_1,x_2,x_3$ by
\[\begin{array}{ccccl}
x_1 &:=& \displaystyle \frac{2K(k)}{\pi} &=& 3,264710703\dots \,,\\ \\
x_2 &:=& \displaystyle \frac{2E(k)}{\pi} &=& 0,637448893\dots \,,\\\\
x_3 &:=& k &=& 0,999718575\dots
\end{array} \]
It is known that the three parameters $x_1,x_2,x_3$ are algebraically independent over ${\Q}$. Furthermore, 
there are two explicitly given polynomials $f_i(X_1,X_2,X_3,Y_i)\in {\Q}[X_1,X_2,X_3,Y_i]$ such that 
$f_i(x_1,x_2,x_3,y_i)=0$ for $i=1,2$. We apply Theorem\,4 with $n=3$, $m=2$ and obtain that $\zeta_{Fib}(4)$ 
and $\zeta_{Fib}(8)$ are algebraically independent over the field ${\Q}(k)$.  
\end{itemize}


\begin{thebibliography}{99}
\bibitem{Baker} A.\,Baker, \textit{Transcendental Number Theory}, Cambridge University Press (1975).
\bibitem{Elsner1} C.\,Elsner, Sh.\,Shimomura, I.\,Shiokawa, and Y.\,Tachiya, \textit{Algebraic independence results for the
sixteen families of \(q\)-series}, Ramanujan Journal {\bf 22} no.\,3, (2010), 315-344.  
\bibitem{Elsner2} C.\,Elsner, Sh.\,Shimomura, I.\,Shiokawa, \textit{A remark on Nesterenko's theorem for Ramanujan functions}, 
Ramanujan Journal {\bf 21} no.\,2, (2010), 211--221.
\bibitem{Elsner3} C.\,Elsner, Sh.\,Shimomura, I.\,Shiokawa, \textit{Algebraic independence results for reciprocal sums of
Fibonacci numbers}, Acta Arithmetica {\bf 148.3} (2011), 205--223.
\bibitem{Elsner5} C.\,Elsner, Sh.\,Shimomura, I.\,Shiokawa, \textit{Algebraic independence of certain numbers related to
modular functions}, Functiones et Approximatio {\bf 47.1} (2012), 121--141; DOI: 10.7169/facm/2012.47.1.10.
\bibitem{Elsner6} C.\,Elsner, \textit{Algebraic independence results for values of theta-constants}, Functiones et Approximatio 
Commentarii Mathematici, Funct. Approx. Comment. Math. {\bf 52(1)} (2015), 7-27.
\bibitem{Elsner8} C.\,Elsner, Yu.V.\,Nesterenko, I.\,Shiokawa, \textit{Algebraic independence of values of exponential type
power series}, Moscow Journal of Combinatorics and Number Theory {\bf 3}, no.\,2-3 (2013).
\bibitem{Elsner10} C.\,Elsner and T.Komatsu,\,{\em A recurrence formula for leaping convergents of non-regular continued fractions},\/
Lin. Alg. Appl, 428 (2008), 824-833.
\bibitem{Hardy} G.H.\,Hardy and E.M.\,Wright, {\em An introduction to the theory of numbers},\/ fifth edition, Clarendon Press, Oxford, 1979.
\bibitem{Lang} S.\,Lang, {\em Algebra}, 3rd ed., Addison-Wesley (1993).
\bibitem{Nesterenko1} Yu.V.\,Nesterenko, \textit{Modular functions and transcendence problems}, Math.\,Sbornik {\bf 187} 
(1996), 65-96; C.R.\,Acad.\,Sci.\,Paris; Ser.\,1 {\bf 322} (1996), 909-914.
\bibitem{Ramanujan} S.\,Ramanujan, {\em On certain arithmetical functions\/}, Trans.\,Cambridge\,Philos.\,Soc. {\bf 22} (1916), 159-184.
\bibitem{Shidlovskii} A.B.\,Shidlovskii, {\em Transcendental Numbers}, Walter de Gruyter (1989).
\bibitem{Stein} M.\,Stein, {\em Algebraic independence results for reciprocal sums of Fibonacci and Lucas numbers},
Dissertation an der Leibniz Universit\a t Hannover (2012), Katalognummer der Techn. Informationsbibliothek Hannover:\,H\,12\,B\,58\,a,b;\,http://edok01.tib.uni-hannover.de/edoks/e01dh12/684662426.pdf
\bibitem{Zariski} O.\,Zariski, P.\,Samuel, {\em Commutative Algebra}, vol.\,1, Graduate Texts in Mathematics~28, Springer (1975).
\end{thebibliography}
\end{document}